\documentstyle{article}
\def\c{\mbox{\sf c}}
\def\e{\mbox{\sf e}}
\def\g{\mbox{\sf g}}
\def\k{\mbox{\sf k}}
\def\x{\mbox{\sf x}}
\def\Q{\mbox{\sf Q}}
\def\R{\mbox{\sf R}}
\begin{document}
\renewcommand{\thefootnote}{\fnsymbol{footnote}}

\begin{center}{\Large\bf q-combinatorics and
quantum integrability}
\bigskip

{\large\bf A. Yu. Volkov\footnote{On leave of absence from
Steklov Mathematical Institute, St. Petersburg}}

{\em ENSLAPP, ENSLyon}
\end{center}
\bigskip
\begin{quote}{\bf ABSTRACT.} The idea that a Dynkin diagram
can provide one of the `spatial' variables for an integrable
difference-difference system is no news. I propose
a `model' where the only variable is of this sort.
\end{quote}
\bigskip

\noindent It has been observed by S.~Fomin and An.~Kirillov
[FK] that if in some algebra there were a bunch of
elements $\R_n(\lambda)$
\begin{description}
\item[(i)] obeying Artin-Yang-Baxter's commutation relations
$$ \R_{n+1}(\lambda-\mu)\R_{n}(\lambda)\R_{n+1}(\mu)
   =\R_{n}(\mu)\R_{n+1}(\lambda)\R_{n}(\lambda-\mu)        $$
$$ \R_m\R_n=\R_n\R_m \qquad
   \mbox{if}\qquad |m-n|\neq 1                             $$
\item[(ii)] and depending on the spectral parameter in an
exponential way
$$ \R_{n}(\lambda)=\R_{n}(\lambda-\mu)\R_{n}(\mu)          $$
\end{description}
then ordered products of those R-matrices would all commute
with each other:
$$ \bigg(\R_N(\lambda)\ldots\R_2(\lambda)\R_1(\lambda)\bigg)
   \bigg(\R_N(\mu)\ldots\R_2(\mu)\R_1(\mu)\bigg)
   \qquad\qquad\qquad                                      $$
$$ \qquad\qquad\qquad =
   \bigg(\R_N(\mu)\ldots\R_2(\mu)\R_1(\mu)\bigg)
   \bigg(\R_N(\lambda)
   \ldots\R_2(\lambda)\R_1(\lambda)\bigg)       .          $$
Unfortunately, this was not to become an algebraic skeleton
of quantum integrability. Conditions (i) and (ii) proved
too restrictive and all the examples that emerged were
of a nondeformable classical sort. Fortunately, it is
reparable. Let us go back to basics and make R-matrices
depend on two spectral parameters so that (i) and (ii) rather
read
$$ \begin{array}{c}
   \R_{n+1}(\lambda,\mu)\R_{n}(\lambda,\nu)\R_{n+1}(\mu,\nu)
   =\R_{n}(\mu,\nu)\R_{n+1}(\lambda,\nu)\R_{n}(\lambda,\mu)\\
   \\ \R_m\R_n=\R_n\R_m \qquad
   \mbox{if}\qquad |m-n|\neq 1
   \end{array}\eqno\mbox{(i)}                              $$
$$ \R_{n}(\lambda,\nu)=\R_{n}(\lambda,\mu)\R_{n}(\mu,\nu).
   \eqno\mbox{(ii)}                                        $$
{\bf PROPOSITION} remains: ordered products
$$ \Q(\lambda,\mu)\equiv\R_N(\lambda,\mu)\ldots
   \R_2(\lambda,\mu)\R_1(\lambda,\mu)                      $$
commute whenever their second `arguments' coincide:
$$ \Q(\lambda,\nu)\Q(\mu,\nu)=\Q(\mu,\nu)\Q(\lambda,\nu).  $$
$$ \begin{array}{l}Proof:\qquad\Q(\lambda,\nu)\Q(\mu,\nu)\\
   \\=\R_N(\lambda,\nu)
   \left(\R_{N-1}(\lambda,\nu)\R_N(\mu,\nu)\right)\ldots
   \left(\R_1(\lambda,\nu)\R_2(\mu,\nu)\right)\R_1(\mu,\nu)\\
   \\=\R_N(\mu,\nu)\R_N(\lambda,\mu)
   \left(\R_{N-1}(\lambda,\nu)\R_N(\mu,\nu)\right)\ldots
   \left(\R_1(\lambda,\nu)\R_2(\mu,\nu)\right)\R_1(\mu,\nu)\\
   \\=\R_N(\mu,\nu)
   \left(\R_{N-1}(\mu,\nu)\R_N(\lambda,\nu)\right)\,\ldots\,
   \left(\R_1(\mu,\nu)\R_2(\lambda,\nu)\right)
   \R_1(\lambda,\mu)\R_1(\mu,\nu)\\
   \\=\R_N(\mu,\nu)
   \left(\R_{N-1}(\mu,\nu)\R_N(\lambda,\nu)\right)\ldots
   \left(\R_1(\mu,\nu)\R_2(\lambda,\nu)\right)
   \R_1(\lambda,\nu)\\
   \\=\Q(\mu,\nu)\Q(\lambda,\nu)
   \qquad\Box        \end{array}                           $$
It is plain to see that this would turn right back into
Fomin-Kirillov's case if I added the
usual (iii) $\R(\lambda,\mu)=\R(\lambda-\mu)$. Naturally,
I do not.
{\bf PROPOSITION}: in the algebra whose only two generators
are bound by Serre-style commutation relations
$$ \x_1\x_1\x_2+\x_2\x_1\x_1 q=\x_1\x_2\x_1(1+q)           $$
$$ \x_1\x_2\x_2+\x_2\x_2\x_1 q=\x_2\x_1\x_2(1+q)           $$
criteria (i) and (ii) are met by the elements
$$ \R_n(\lambda,\mu)=(\x_n)^\mu_\lambda
   \equiv\prod^{\lambda-1}_{j=\mu}(1-\x_n q^j).            $$
(ii) comes free of charge, {\em proof} of (i) (only the first
line applies) starts with {\bf lemma} establishing something
like Campbell-Hausdorff multiplication rules:
$$ (\x_1)^\mu_\lambda(\x_2)^\mu_\lambda
   =\prod^{\lambda-1}_{j=\mu}
   \left(1+\c q^{2j}-(\x_1+\x_2+\k q^\lambda) q^j\right)   $$
$$ (\x_2)^\mu_\lambda(\x_1)^\mu_\lambda
   =\prod^{\lambda-1}_{j=\mu}
   \left(1+\c q^{2j}-(\x_1+\x_2+\k q^\mu) q^j\right)       $$
with
$$ \k=\frac{\x_1\x_2-\x_2\x_1}{1-q}\qquad\qquad
   \c=\frac{\x_1\x_2-\x_2\x_1 q}{1-q},                     $$
the element $\c$ actually being central. Proof is by
induction for there is nothing but polynomials in here.
Let me omit it. So,
$$ \begin{array}{l}
   (\x_2)^\mu_\lambda(\x_1)^\nu_\lambda(\x_2)^\nu_\mu
   =(\x_2)^\mu_\lambda(\x_1)^\mu_\lambda
   (\x_1)^\nu_\mu(\x_2)^\nu_\mu\\\\
   =\prod^{\lambda-1}_{j=\nu}
   \left(1+\c q^{2j}-(\x_1+\x_2+\k q^\mu) q^j\right)\\\\
   =(\x_1)^\nu_\mu(\x_2)^\nu_\mu
   (\x_2)^\mu_\lambda(\x_1)^\mu_\lambda
   =(\x_1)^\nu_\mu(\x_2)^\nu_\lambda(\x_1)^\mu_\lambda
   \qquad\Box \end{array}                                  $$
The two propositions combine into the message of this note:
the algebra whose $r$ generators commute like this
$$ \begin{array}{l}\x_n\x_n\x_{n+1}+\x_{n+1}\x_n\x_n q
   =\x_n\x_{n+1}\x_n(1+q)\\ \\
   \x_n\x_{n+1}\x_{n+1}+\x_{n+1}\x_{n+1}\x_n q
   =\x_{n+1}\x_n\x_{n+1}(1+q)\\ \\
   \x_m\x_n=\x_n\x_m \qquad
   \mbox{if}\qquad |m-n|\neq 1  \end{array}                $$
contains a good supply
$$ \Q(\lambda)=(\x_r)^{}_\lambda\ldots
   (\x_2)^{}_\lambda(\x_1)^{}_\lambda\qquad\qquad
   (\cdot)^{}_\lambda\equiv(\cdot)^0_\lambda               $$
of mutually commuting elements
$$ \Q(\lambda)\Q(\mu)=\Q(\mu)\Q(\lambda).                  $$
In conclusion, some remarks. The definition of
$(\cdot)^\mu_\lambda$ required integer $\lambda,\mu$ bound
by $\lambda\geq\mu$. It would be more
practical (if less stylish) to do without those limitations.
Formal power series in $\x$'s could help. An obvious
identity
$$ (\x)^\mu_\lambda=\frac{(\x q^\mu)^{}_\infty}
   {(\x q^\lambda)^{}_\infty}\qquad\qquad
   (\x)^{}_\infty\equiv
   \prod^{\infty}_{j=0}(1-\x q^j)                          $$
would then double as a definition of $(\cdot)^\mu_\lambda$
for non-integer $\lambda,\mu$. As a matter of fact, AYB
relations
$$ (\x_{n+1})^\mu_\lambda(\x_n)^\nu_\lambda
   (\x_{n+1})^\nu_\mu=(\x_n)^\nu_\mu
   (\x_{n+1})^\nu_\lambda(\x_n)^\mu_\lambda                $$
survive the extrapolation. Anyway, it is perhaps more
important to guess where the commutation relations
governing the $\x$'s belong. Quantum $A_r$ algebra provides
$$ \begin{array}{l}\e_n\e_n\e_{n+1}+\e_{n+1}\e_n\e_n
   =\e_n\e_{n+1}\e_n(q^{\frac{1}{2}}+q^{-\frac{1}{2}})\\ \\
   \e_n\e_{n+1}\e_{n+1}+\e_{n+1}\e_{n+1}\e_n
   =\e_{n+1}\e_n\e_{n+1}(q^{\frac{1}{2}}+q^{-\frac{1}{2}})
                                \end{array}                $$
and with a little help of suitable `quantum coordinates'
$$ \chi_m\e_n=\e_n\chi_m\qquad\qquad
   \chi_n\chi_{n+1}=q^{\frac{1}{2}}\chi_{n+1}\chi_n
   \qquad\qquad n=1,2,\ldots,r-1                           $$
`vectors' $\chi_n\e_n$ just fit the commutation relations
prescribed for the $\x$'s. $\Q(\infty)$ becomes a piece
of `quantized' and `bosonized' Gauss decomposition
$$ \begin{array}{rrllll}\g&=&(\chi_r\e_r)^{}_\infty&\ldots&
   (\chi_2\e_2)^{}_\infty&(\chi_1\e_1)^{}_\infty\\\\
   &\times&(\chi_{2r-1}\e_r)^{}_\infty&\ldots&
   (\chi_{r+1}\e_2)^{}_\infty&\\
   &&\vdots&&&\\
   &\times&(\chi_{r(r+1)/2}\e_r)^{}_\infty&&&\\\\
   &\times&\g_{d}\g_{l}&&&\end{array}                      $$
a l\`a Morozov-Vinet [MV]. Let me decline
further comments on this issue.
Instead, let me mention that the
same `R-matrix' $(\cdot)^\mu_\lambda$ used along the
guidelines of [V] provides the quantization of a very
major nonlinear difference-difference system [H]
$$ 1+\psi_{m+1,n}\psi_{m,n+1}+\psi_{m+1,n+1}\psi_{m,n}=-
   \Lambda\psi_{m+1,n}\psi_{m,n+1}\psi_{m+1,n+1}\psi_{m,n}.$$
The case $\Lambda=1$ is actually
linear, $\Lambda=0$ (corresponding, by the way,
to $\lambda=\infty$) approximates the Liouville equation
while everything else is the sine-Gordon equation. I think
we've got a few more pieces of the big puzzle called
Quantum Solitons.
\bigskip
\begin{quote}
{\bf Acknowledgements.}
I would like to thank L.~Faddeev, An.~Kirillov,
R.~Kashaev, J.-M.~Maillet, A.~Morozov and V.~Tarasov
for stimulating discussions.
\end{quote}
\section*{\normalsize\bf References}

\begin{itemize}
  \item[{[FK]}] S.~Fomin and An.~Kirillov, Discrete
Mathematics 153 (1996).
  \item[{[MV]}] A.~Morozov and L.~Vinet, hep-th/9409093.
  \item[{[V]}] A.~Yu.~Volkov, hep-th/9512024.
  \item[{[H]}] R.~Hirota, J. Phys. Soc. Japan 43 (1977).
\end{itemize}
\end{document}